\documentclass[10pt]{article}
\usepackage{amsmath,amsthm,amsfonts,amssymb,amscd,enumerate,pictex}

\newtheorem{SConjecture}[subsection]{Singer's Conjecture}
\newtheorem{subSCConjecture}[subsection]{Singer's Conjecture for Coxeter groups}

\newtheorem{Mainthm}[subsection]{Main Theorem}
\newtheorem{Theorem}[subsection]{Theorem}

\newtheorem{Cor}[subsection]{Corollary}
\newtheorem{Lemma}[subsection]{Lemma}

\newtheorem{Example}[subsection]{Example}

\DeclareMathOperator{\Vol}{Vol}

\newcommand{\cs}{\mathcal{S}}

\newcommand{\cH}{\mathcal{H}}

\newcommand{\Ltwo}{L^2}
\newcommand{\mfh}{\mathfrak{h}}

\newcommand{\ltwo}{\ell^2}

\newcommand{\St}{\operatorname{St}}

\newcommand{\Edge}{\operatorname{Edge}}

\newcommand{\gS}{\Sigma}

\newcommand{\ga}{\alpha}
\newcommand{\gb}{\beta}

\newcommand{\BS}{\mathbb{S}}
\newcommand{\BH}{\mathbb{H}}
\newcommand{\BZ}{\mathbb{Z}}
\newcommand{\BR}{\mathbb{R}}
\newcommand{\BE}{\mathbb{E}}
\newcommand{\BN}{\mathbb{N}}

\newenvironment{enumeratei}{\begin{enumerate}[\upshape (i)]}
        {\end{enumerate}}
\newenvironment{enumeratea}{\begin{enumerate}[\upshape 
(a)]}{\end{enumerate}}

\numberwithin{equation}{section}
\begin{document}

\title{Geometrization of 3-dimensional Coxeter orbifolds and Singer's conjecture}

\author{Timothy A. Schroeder}

\date{\today}
\maketitle

\begin{abstract}
Associated to any Coxeter system $(W,S)$, there is a labeled simplicial complex $L$ and a contractible CW-complex $\gS_L$ (the Davis complex) on which $W$ acts properly and cocompactly.  $\gS_L$ admits a cellulation under which the nerve of each vertex is $L$.  It follows that if $L$ is a triangulation of $\BS^{n-1}$, then $\gS_L$ is a contractible $n$-manifold.  In this case, the orbit space, $K_L:=\gS_L/W$, is a \emph{Coxeter orbifold}.  We prove a result analogous to the JSJ-decomposition for $3$-dimensional manifolds: Every $3$-dimensional Coxeter orbifold splits along Euclidean suborbifolds into the \emph{characteristic suborbifold} and simple (hyperbolic) pieces.  It follows that every $3$-dimensional Coxeter orbifold has a decomposition into pieces which have hyperbolic, Euclidean, or the geometry of $\BH^2\times\BR$.  (We leave out the case of spherical Coxeter orbifolds.)  A version of Singer's conjecture in dimension $3$ follows: That the reduced $\ltwo$-homology of $\gS_L$ vanishes.  
\end{abstract}

\section{Introduction}
\label{intro}
The following conjecture is attributed to Singer.
\begin{SConjecture}\label{conj:singer} If $M^{n}$ is a closed aspherical manifold, then the reduced $\ltwo$-homology of $\widetilde{M}^n$, $\cH_{i}(\widetilde{M}^{n})$, vanishes for all $i\neq\frac{n}{2}$.
\end{SConjecture}
For details on $\ltwo$-homology theory, see \cite{davismoussong}, \cite{do2} and \cite{eckmann}, which is particularly useful and easy to read.  

Let $S$ be a finite set of generators.  A \emph{Coxeter matrix} on $S$ is a symmetric $S\times S$ matrix $M=(m_{st})$ with entries in $\BN\cup\{\infty\}$ such that each diagonal entry is $1$ and each off diagonal entry is $\geq 2$.  The matrix $M$ gives a presentation of an associated \emph{Coxeter group} $W$:
\begin{equation}\label{e:coxetergroup}
	W=\left\langle S\mid (st)^{m_{st}}=1, \text{ for each pair } (s,t) \text{ with } m_{st}\neq\infty\right\rangle.
\end{equation}
The pair $(W,S)$ is called a \emph{Coxeter system}.  Denote by $L$ the nerve of $(W,S)$.  ($L$ is a simplicial complex with vertex set $S$, the precise definition will be given in section \ref{ss:davis}.)  In several papers (e.g., \cite{davisannals}, \cite{davisbook}, and \cite{davismoussong}), M. Davis describes a construction which associates to any Coxeter system $(W,S)$, a simplicial complex $\gS(W,S)$, or simply $\gS$ when the Coxeter system is clear, on which $W$ acts properly and cocompactly.  The two salient features of $\gS$ are that (1) it is contractible and (2) that it admits a cellulation under which the nerve of each vertex is $L$.  It follows that if $L$ is a triangulation of $\BS^{n-1}$, $\gS$ is an $n$-manifold.  

If $G$ is a torsion-free subgroup of finite index in $W$, then $G$ acts freely on $\gS$ and $\gS/G$ is a finite complex.  By $(1)$, $\gS/G$ is aspherical.  Hence, if $L$ is homeomorphic to an $(n-1)$-sphere, Davis' construction gives examples of closed aspherical $n$-manifolds and Conjecture \ref{conj:singer} for such manifolds becomes the following.

\begin{subSCConjecture}\label{conj:singerc} Let $(W,S)$ be a Coxeter group such that its nerve, $L$, is a triangulation of $\BS^{n-1}$.  Then $\cH_{i}(\gS_L)=0$ for all $i\neq\frac{n}{2}$.
\end{subSCConjecture}

Conjecture \ref{conj:singer} holds for elementary reasons in dimensions $\leq 2$.  In \cite{LL}, Lott and L\"uck prove that Conjecture \ref{conj:singer} holds for those aspherical $3$-manifolds for which Thurston's Geometrization Conjecture is true.  (Hence, by Perelman, all aspherical $3$-manifolds.)  Thurston proved in \cite{thurston} that the Geometrization Conjecture holds for Haken $3$-manifolds; and in \cite{do2}, Davis and Okun show that when $(W,S)$ is \emph{right-angled} (this means that generators either commute, or have no relation), Davis' construction yields examples of Haken $3$-manifolds.  Thus, they show that Thurston's Geometrization Conjecture holds for closed aspherical $3$-manifolds arising as quotient spaces of right-angled Davis complexes.  Also in \cite{do2}, the authors show that if Conjecture \ref{conj:singerc} for right-angled Coxeter systems is true in some odd dimension $n$, then it is also true in dimension $n+1$.  Hence, the Lott and L\"uck result implies that Conjecture \ref{conj:singerc} for right-angled Coxeter systems is true for $n=3$ and, therefore, also for $n=4$.  

In the case $L$ is a triangulation of $\BS^{n-1}$, $K_L:=\gS_L/W$ is an $n$-dimensional \emph{Coxeter orbifold}.  A Coxeter orbifold is an orbifold with underlying space an $n$-disk such that each local isotropy group is a finite reflection group.  We call a Coxeter orbifold \emph{spherical, Euclidean, or hyperbolic} if it is the quotient of a reflection group on the sphere, Euclidean space or hyperbolic space, respectively.  In this paper, we deal exclusively with non-spherical Coxeter orbifolds.  

In \cite[Theorem 2]{andreev2}, E. Andreev gives the necessary and sufficient conditions for an abstract $3$-dimensional polytope $P$ with assigned dihedral angles to be realized as a convex polytope in $\BH^3$.  Andreev's theorem guides us to our main result: A decomposition of Coxeter orbifolds analogous to the JSJ-decomposition for $3$-manifolds (W. Jaco and P. Shalen \cite{jash}, K. Johannson \cite{joh}). 
\begin{Mainthm}\label{t:main} Let $K$ be a closed, irreducible $3$-dimensional Coxeter orbifold.  There is a unique minimal collection of Euclidean suborbifolds of $K$ forming the boundary of the characteristic suborbifold of $K$.  
\end{Mainthm}
We then show that Conjecture \ref{conj:singerc} in dimension $3$ follows from Theorem \ref{t:main}.  For further information on $3$-dimensional orbifolds, the reader is directed to \cite[Chapter 13: Orbifolds]{thurston2}.

\section{The Davis complex and Coxeter orbifolds}\label{s:sec1}
Let $(W,S)$ be a Coxeter system.  Given a subset $U$ of $S$, define $W_{U}$ to be the subgroup of $W$ generated by the elements of $U$.  A subset $T$ of $S$ is \textit{spherical} if $W_T$ is a finite subgroup of $W$.  In this case, we will also say that the subgroup $W_{T}$ is spherical.  Denote by $\cs$ the poset of spherical subsets of $S$, partially ordered by inclusion.  Given a subset $V$ of $S$, let $\cs_{\geq V}:=\{T\in \cs|V\subseteq T\}$.  Similar definitions exist for $<, >, \leq$.  For any $w\in W$ and $T\in \cs$, we call the coset $wW_{T}$ a \emph{spherical coset}.  The poset of all spherical cosets we will denote by $W\cs$.  

\subsection{The Davis complex}\label{ss:davis}
Let $K=|\cs|$, the geometric realization of the poset $\cs$.  It is a finite simplicial complex.  Denote by $\gS(W,S)$, or simply $\gS$ when the system is clear, the geometric realization of the poset $W\cs$.  This is the Davis complex.  The natural action of $W$ on $W\cs$ induces a simplicial action of $W$ on $\gS$ which is proper and cocompact.  $\gS$ is a model for $\underline{E}W$, a \emph{universal space for proper $W$-actions}.  (See Definition \cite[2.3.1]{davisbook}.)  $K$ includes naturally into $\gS$ via the map induced by $T\rightarrow W_{T}$, so we view $K$ as a subcomplex of $\gS$.  Note that $K$ is a strict fundamental domain for the action of $W$ on $\gS$.  

The poset $\cs_{>\emptyset}$ is an abstract simplicial complex.  This simply means that if $T\in\cs_{>\emptyset}$ and $T'$ is a nonempty subset of $T$, then $T'\in \cs_{>\emptyset}$.  Denote this simplicial complex by $L$, and call it the \emph{nerve} of $(W,S)$.  The vertex set of $L$ is $S$ and a non-empty subset of vertices $T$ spans a simplex of $L$ if and only if $T$ is spherical.  Define a labeling on the edges of $L$ by the map $m:\Edge(L)\rightarrow \{2,3,\ldots\}$, where $\{s,t\}\mapsto m_{st}$.  This labeling accomplishes two things: (1) the Coxeter system $(W,S)$ can be recovered (up to isomorphism) from $L$ and (2) the $1$-skeleton of $L$ inherits a natural piecewise spherical structure in which the edge $\{s,t\}$ has length $\pi-\pi/m_{st}$.  $L$ is then a \emph{metric flag} simplicial complex (see Definition \cite[I.7.1]{davisbook}).  This means that any finite set of vertices, which are pairwise connected by edges, spans a simplex of $L$ if an only if it is possible to find some spherical simplex with the given edge lengths.  In other words, $L$ is ``metrically determined by its $1$-skeleton.''  

For the purpose of this paper, we will say that labeled (with integers $\geq 2$) simplicial complexes are \emph{metric flag} if they correspond to the labeled nerve of some Coxeter system.  We will often indicate these complexes simply with their $1$-skeleton, understanding the underlying Coxeter system and Davis complex.  We write $\gS_L$ to denote the Davis complex associated to the nerve $L$ of $(W,S)$.  

\paragraph{A cellulation of $\gS$ by Coxeter cells.}  $\gS$ has a coarser cell structure: its cellulation by ``Coxeter cells.''  (References for this section include \cite{davisbook} and \cite{do2}.)  The features of the Coxeter cellulation are summarized by \cite[Proposition 7.3.4]{davisbook}.  We note here that, under this cellulation, the link of each vertex is $L$.  It follows that if $L$ is a triangulation of $\BS^{n-1}$, then $\gS$ is a topological $n$-manifold.  

\paragraph{Special subcomplexes.}  Suppose $A$ is a full subcomplex of $L$.  Then $A$ is the nerve for the subgroup generated by the vertex set of $A$.  We will denote this subgroup by $W_{A}$.  (This notation is natural since the vertex set of $A$ corresponds to a subset of the generating set $S$.)  Let $\cs_{A}$ denote the poset of the spherical subsets of $W_A$ and let $\gS_{A}$ denote the Davis complex associated to $(W_{A},A^{(0)})$.  The inclusion $W_{A}\hookrightarrow W_{L}$ induces an inclusion of posets $W_{A}\cs_{A}\hookrightarrow W_{L}\cs_{L}$ and thus an inclusion of $\gS_{A}$ as a subcomplex of $\gS_{L}$.  Such a subcomplex will be called a \emph{special subcomplex} of $\gS_{L}$.  Note that $W_{A}$ acts on $\gS_{A}$ and that if $w\in W_{L}-W_{A}$, then $\gS_{A}$ and $w\gS_{A}$ are disjoint copies of $\gS_{A}$.  Denote by $W_{L}\gS_{A}$ the union of all translates of $\gS_{A}$ in $\gS_{L}$.  

\paragraph{A mirror structure on $K$.} If $L$ is the triangulation of an $(n-1)$-sphere, then we have a another cellulation of $K$ and $\gS$.  For each $T\in\cs$, let $K_T$ denote the geometric realization of the subposet $\cs_{\geq T}$.  $K_T$ is a triangulation of a $k$-cell, where $k=n-|T|$.  We then define a new cell structure on $K$ by declaring the family $\{K_T\}_{T\in\cs}$ to be the set of cells in $K$.  Under the $W$-action on $\gS$, the finite subgroup $W_T$ is the stabilizer of the cell $K_T$.  We write $K_L$ to indicate $K$ equipped with this cellulation and note that it extends to a cellulation of $\gS_L$.  $K_L$ is an $n$-dimensional \emph{Coxeter orbifold}: $K_L$ is an orbifold with underlying space an $n$-disk such that isotropy group of every face is a finite reflection group.  

\subsection{Coxeter orbifolds}  A $3$-dimensional Coxeter orbifold $K$ is said to be \emph{irreducible} if every $2$-dimensional, spherical Coxeter suborbifold bounds the quotient of a $3$-disk by a finite reflection group.  Every $3$-dimensional $K$ can be decomposed along spherical suborbifolds into irreducible pieces (i.e. $K$ is a connected sum of irreducible Coxeter orbifolds).  We say a face of $K$ is \emph{labeled} with $G$ if $G$ is the isotropy group of this face.  $K$ is \emph{closed} if every codimension one face of $K$ is labeled with $\BZ_2$.  A Coxeter orbifold $K$ is \emph{atoroidal} if it has no incompressible, $2$-dimensional Euclidean Coxeter suborbifolds.  The \emph{characteristic suborbifold} of $K$ is the minimal (possibly disconnected) suborbifold containing all Euclidean suborbifolds, i.e. its complement is atoroidal.  

\paragraph{Andreev's Theorem.}  In \cite{andreev2}, Andreev gives the necessary and sufficient conditions for abstract $3$-dimensional polytopes, with assigned dihedral angles in $\left(0,\frac{\pi}{2}\right]$, to be realized as (possibly ideal) convex polytopes in $\BH^3$ (these conditions are listed below, Theorem \ref{t:2:andreev}).  In order for this convex polytope to tile $\BH^3$, the assigned dihedral angles must be integer submultiples of $\pi$.  

As described in Section \ref{ss:davis}, any metric flag triangulation $L$ of $\BS^2$ defines a closed, irreducible $3$-dimensional Coxeter orbifold $K_L$.  The boundary complex of $K_L$ is combinatorially dual to $L$, so $K_L$ has codimension $1$ faces corresponding the elements of $S$.  In fact, if $Z$ is \emph{any} (labeled) cell complex homeomorphic to $\BS^2$, in the strict sense that any non-empty intersection of two cells is a cell, then $Z$ is combinatorially dual to the boundary complex of a $3$-dimensional convex polytope, which we will denote by $K_{Z}$.  (In this generality, $K_Z$ is not necessarily a Coxeter orbifold.  The subscript indicates the dual cellulation of $\BS^2$.)  We assign dihedral angles to $K_Z$ so that the angle between faces dual to vertices $s$ and $t$ is $\pi/m_{st}$, where $m_{st}$ is the label on the edge between $s$ and $t$.   
\begin{Theorem}\label{t:2:andreev} \textup{(\cite[Theorem 2]{andreev2})} Let $P$ be an abstract three-dimensional polyhedron, not a simplex, such that three or four faces meet at every vertex.  The following conditions are necessary and sufficient for the existence in $\BH^{3}$ of a convex polytope of finite volume of the combinatorial type $P$ with the dihedral angles $\ga_{ij}\leq\frac{\pi}{2}$ (where $\ga_{ij}$ is the dihedral angle between the faces $F_{i},F_{j}$):
\begin{enumeratei}
	\item\label{i:i} If $F_{1}, F_{2}$ and $F_{3}$ are all the faces meeting at a vertex of $P$, then $\ga_{12}+\ga_{23}+\ga_{31}\geq \pi$; and if $F_{1}, F_{2}, F_{3}, F_{4}$ are all the faces meeting at a vertex of $P$ then $\ga_{12}+\ga_{23}+\ga_{34}+\ga_{41}=2\pi$.
	\item\label{i:ii} If three faces intersect pairwise but do not have a common vertex, then the angles at the three edges of intersection satisfy $\ga_{12}+\ga_{23}+\ga_{31}<\pi$.
	\item\label{i:iii} Four faces cannot intersect cyclically with all four angles $=\pi/2$ unless two of the opposite faces also intersect.  
	\item\label{i:iv} If $P$ is a triangular prism, then the angles along the base and top cannot all be $\frac{\pi}{2}$.  
	\item\label{i:v} If among the faces $F_{1},F_{2},F_{3}$ we have $F_{1}$ and $F_{2}$, $F_{2}$ and $F_{3}$ adjacent, but $F_{1}$ and $F_{3}$ not adjacent, but concurrent at one vertex and all three do not meet in one vertex, then $\ga_{12}+\ga_{23}<\pi$.
\end{enumeratei}
\end{Theorem} 

Now, with $L$ a metric flag triangulation of $\BS^2$, the conditions of Andreev's Theorem refer to certain configurations of $L$ which, in turn, correspond to certain suborbifolds of $K_L$.  In the next section, we'll identify which subcomplexes of $L$ define components of the characteristic suborbifold of $K_L$.
 

\section{The Geometrization of $K_L$}\label{s:geom}
Let $L$ be a metric flag triangulation of $\BS^2$, and unless otherwise noted, not the boundary of a $3$-simplex.  Let $K_L$ denote the corresponding Coxeter orbifold.

\paragraph{Euclidean vertices.}  If $s$ is a vertex of $L$, define the \emph{link of $s$ in $L$}, $L_{s}$, to be the subcomplex of $L$ consisting of all closed simplices which are contained in simplices containing $s$, but do not themselves contain $s$.   Define the \emph{star of $s$ in $L$}, $\St_{L}(s)$, to be the subcomplex of $L$ consisting of all closed simplices which contain $s$.  

The \emph{valence} of a vertex $s$ is the number of vertices in its link.  We say that a vertex $s$ is \emph{3-Euclidean} if $s$ has valence $3$ and if $s_0, s_1, s_2$ are the vertices in this link, then 
\[\frac{\pi}{m_{s_0s_1}}+\frac{\pi}{m_{s_1s_2}}+\frac{\pi}{m_{s_2s_0}}=\pi.\]
We say that $s\in T$ is \emph{4-Euclidean}, if $s$ has valence $4$ and if $s_{0}, s_{1}, s_{2}, s_{3}$ are the vertices in this link, then $m_{s_{i}s_{i+1}}=2$ for $i=0,1,2,3$ (mod($4$)).  We'll say that the vertex $s$ is \emph{Euclidean} if it is either $3$- or $4$-Euclidean.  

\begin{Lemma}\label{l:full} Let $s$ be a Euclidean vertex.  
\begin{enumeratea}
	\item If $s$ is a $3$-Euclidean vertex, then $L_s$ and $\St_L(s)$ are full subcomplexes of $L$.
	\item If $s$ is a $4$-Euclidean vertex and $L$ is not the suspension of a $3$-gon, then $L_s$ and $\St_L(s)$ are full subcomplexes of $L$.
\end{enumeratea}  
\end{Lemma}
\begin{proof} $(a)$: This is immediate since $L$ is not the boundary of a $3$-simplex.\\
$(b)$: For a $4$-Euclidean vertex $s$, $L_s$ and $\St_L(s)$ can only fail to be full if $L$ is the suspension of a $3$-gon.
\end{proof}

\paragraph{The geometry of the stars of Euclidean vertices.}  Suppose $s$ is a Euclidean vertex of $L$ (if $s$ is $4$-Euclidean, require that $L$ is not the suspension of a $3$-gon).  If each edge in $(\St_{L}(s)-L_{s})$ is labeled $2$, then $W_{\St_L(s)}=W_{L_s}\times W_s$ and $\gS_{\St_{L}(s)}=\BR^2\times\left[-1,1\right]$ ($\gS_{s}=\left[-1,1\right]$).  We refer to such stars as \emph{right-angled cones}, or \emph{RA-cones}.  The corresponding suborbifold of $K_L$ is a triangular or rectangular prism with one base labeled with the trivial group.  It is the quotient of the $W_{\St_L(s)}$ action on $\gS_{\St_L(s)}$.  Otherwise, let $\left[\St\right]$ denote the complex obtained by capping off the boundary of $\St_{L}(s)$ with a triangular or square cell.  

If $s$ is $3$-Euclidean and the edges in $(\St_{L}(s)-L_s)$ are not all labeled $2$, then the reflection group $W_{\St_{L}(s)}$ is one of the Coxeter groups shown in Figure 6.3 of \cite{humphreys}, the non-compact hyperbolic Coxeter groups ($n=4$).  It acts properly as a classical reflection group on $\BH^3$ with fundamental chamber $K_{\left[\St\right]}$, a simplex of finite volume with one ideal vertex corresponding to the added triangular face of $\left[\St\right]$.  The corresponding suborbifold of $K_L$ is obtained by cutting off the ideal vertex and labeling the resulting face with the trivial group.  

If $s$ is $4$-Euclidean, then the only condition of Theorem \ref{t:2:andreev} $K_{\left[St\right]}$ may fail to satisfy is (\ref{i:v}).  If $K_{\left[St\right]}$ satisfies this condition, then $W_{\St_L(s)}$ acts on $\BH^3$ with finite volume fundamental chamber $K_{\left[St\right]}$.  Again, the corresponding suborbifold is formed by cutting off the ideal vertex and labeling the resulting face with the trivial group.  Otherwise, $\St_L(s)$ is an ``infinite right-angled suspension,'' a case we describe below.

\paragraph{Euclidean circuits.}  Let $C$ be a $3$-circuit in $L$ and let $s_0,s_1,s_2$ be the vertices in this circuit.  We say that $C$ is a \emph{Euclidean $3$-circuit} if  
\[\frac{\pi}{m_{s_0s_1}}+\frac{\pi}{m_{s_1s_2}}+\frac{\pi}{m_{s_2s_0}}=\pi.\]
We say $C$ is an \emph{empty} Euclidean $3$-circuit if $C$ is not the boundary of RA-cone.  It follows from $L$ being metric flag that any Euclidean $3$-circuit $C$ is a full subcomplex.

Let $C$ be a $4$-circuit in $L$.  Order the vertices in this circuit $s_{0},s_{1},s_{2},s_{3}$ so that $s_{i}$ and $s_{i+1}$ are connected by an edge of the circuit and $s_{i}$ and $s_{i+2}$ are not connected by an edge of the circuit ($i=0,1,2,3$ mod($4$)).  We say $C$ is a \emph{Euclidean $4$-circuit} if $m_{s_{i}s_{i+1}}=2$ ($i=0,\ldots,3$ mod($4$)) and if $C$ is not the boundary of the union of two adjacent $2$-simplices $C$.  It follows from $L$ being metric flag that any Euclidean $4$-circuit is a full subcomplex.  Note that by this definition, suspensions of $3$-gons contain no Euclidean $4$-circuits.

For any Euclidean $3$- or $4$-circuit $C$, it is clear $\gS_C=\BE^2$.  These correspond to incompressible, $2$-dimensional Euclidean suborbifolds of $K_L$, so these circuits define part of the characteristic suborbifold of $K_L$.  However, so that this characteristic suborbifold is $3$-dimensional, we require that for $C\subset L$ the corresponding suborbifold of $K_L$ is a triangular or rectangular prism with both bases labeled with the trivial group.  In other words, the suborbifold is the quotient of the $W_C\times\{e\}$ action on $\BE^2\times\left[-1,1\right]$, where $\{e\}$ represents the trivial group.

\paragraph{Right-angled suspensions.}  If a subcomplex $T$ of $L$ is a suspension, i.e. $T=Z\ast P$ where $P$ denotes two points not connected by an edge, and if each suspension edge is labeled $2$, then we call $T$ a \emph{right-angled suspension}, or \emph{RA-suspension}.  $T$ is \emph{maximal} if $Z$ is a full subcomplex of $L$ and if $T$ is not properly contained (as a subcomplex) in another RA-suspension.  $T$ is \emph{infinite} if it is not the suspension of a single edge or of a single vertex.  

\begin{figure}[h]
	\hspace{3.5cm}
	\input{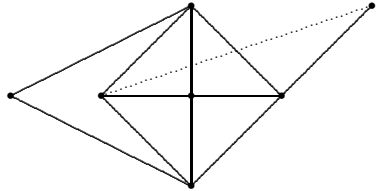}
	\caption{$T$ a RA-cone on a Euclidean $4$-circuit.}
	\label{fig:rtangledcone}
\end{figure}

The maximality condition is well-defined.  Indeed, suppose $T\subsetneq L$ is an infinite RA-suspension that does not have uniquely defined suspension points.  Then since $L$ triangulates $\BS^2$, either $T$ is a RA-cone on a Euclidean $4$-circuit or $T$ is itself a Euclidean $4$-circuit.  Assume $T$ is a RA-cone on a Euclidean $4$-circuit and that in $L$, a point not in $T$ is suspended from each pair of opposite corners of the boundary of $T$.  Since $L$ triangulates $\BS^2$, it must be the same point suspended from each pair.  ($T$ is a disk and each new $4$-circuit bounds a disk in $L$.  The resulting configuration cannot be a subcomplex of $\BS^2$, see Figure \ref{fig:rtangledcone}.)  Therefore, since $L$ is metric flag, $T$ is a ``right-angled octahedron'' (the suspension of a $4$-gon, all labels $2$), and all of $L$.  The geometry of $K_L$ in this case is known: $\gS_L=\BE^3$.  

Now suppose $T$ is a Euclidean $4$-circuit and that in $L$, points not in $T$ are suspended from opposite corners of the $4$-circuit.  If each pair suspends one and the same point, then $T$ is the boundary of a RA-cone on a Euclidean $4$-circuit.  So, assume distinct points are suspended from opposite corners of $T$.  Since $L$ triangulates $\BS^2$, points suspended from different suspension points cannot be connected by an edge.  For then $L$ contains the join of 3 points with 3 points as a subcomplex.  Thus, $T$ and the additional suspended points can be decomposed (uniquely) as the union of two RA-suspensions, glued along $T$.  

\begin{Lemma}\label{l:rasfull} Suppose that $L$ is a metric flag triangulation of $\BS^2$.  Then infinite, maximal RA-suspensions are full subcomplexes of $L$.  
\end{Lemma}
\begin{proof} Suppose that $T$ is an infinite, maximal RA-suspension with suspension points $t$ and $t'$ connected by an edge in $L$.  $T$ contains at least three suspended points, say $s$, $u$ and $v$.  ($T$ is infinite, so it is not the suspension of one point; nor can there be only 2 since then $T$ is a subcomplex of the RA-suspension of the edge connecting $t$ and $t'$, contradicting maximality.)  Then since $L$ is metric flag, three $2$-simplices with vertex sets $\{t,s,t'\}$, $\{t,u,t'\}$ and $\{t,v,t'\}$ are connected along the edge $\{t,t'\}$.  A contradiction to $L$ triangulating $\BS^2$. 
\end{proof}

\paragraph{The geometry of a RA-suspension.}  Let $T=Z\ast P$ be an infinite, maximal RA-suspension, not a Euclidean $4$-circuit nor a RA-cone on a Euclidean $4$-circuit.  Then $W_T=W_Z\times W_P$, where $Z$, a full subcomplex of $L$ by the maximality of $T$, is the nerve of an infinite reflection group.  If $T=L$, then $\gS_L=\BE^3$ or $\BH^2\times\BE$ with compact fundamental chamber the Coxeter orbifold $K_L$.  Otherwise, let $\left[T\right]$ denote the cell complex, homeomorphic to $\BS^2$ obtained by capping the $4$-circuits in the boundary of $T$ with square cells.  Then $W_T$ acts on $\BH^2\times\BE$ with finite volume fundamental chamber $K_{\left[T\right]}$.  Which, when projected to $\BH^2$, has ideal vertices corresponding to the added square faces.  The corresponding suborbifold of $K_L$ is obtained by cutting off the ideal vertices and labeling the resulting faces with the trivial group.  An intermediate cover of this orbifold is isomorphic to $M^2\times \BS^1$, where $M^2$ denotes a surface.  So we say that RA-suspensions are \emph{Seifert fibered.}

\paragraph{Seifert subcomplex.}  We define an equivalence relation on the set of infinite, maximal RA-suspensions in $L$ as follows.  Let $S$ and $S'$ be two such suspensions.  We say $S\sim S'$ if $S\cap S'$ is a Euclidean $4$-circuit.  We call the classes formed by the equivalence relation generated by this symmetric relation the \emph{Seifert subcomplexes} of $L$.  In other words, a Seifert subcomplex is a maximal union of infinite, maximal RA-suspensions glued along $4$-circuits in their boundaries.  The boundary of the resulting Seifert subcomplex is made up of Euclidean $4$-circuits.  (By maximality and by Lemma \ref{l:rasfull}, a $4$-circuit in the boundary cannot bound two simplices.)  Note that RA-suspensions of disjoint points can be Seifert subcomplexes.  Also note that by maximality, a RA-cone on a $4$-Euclidean circuit which is a maximal RA-suspension makes up the entirety of a Seifert subcomplex.  

\paragraph{The geometry of a Seifert subcomplex.}  Let $D$ be a Seifert subcomplex in $L$, where $D$ is not a single RA-suspension.  The geometry of each infinite, maximal RA-suspension $T\subset D$ is $\BH^2\times\BE$.  However, $T$ is glued to another RA-suspension $S$ along a Euclidean $4$-circuit so that suspension points are not identified.  Hence, their corresponding suborbifolds are glued along Euclidean patches with a twist, where the Euclidean factors determined by the suspension points are orthogonal.  In other words, intermediate covers of the orbifolds corresponding to each suspension are glued along tori so that a meridian of one is identified with a longitude of the other.  So the whole of $D$ does not correspond to a Seifert fibered manifold.  However, cutting $D$ along the intersections of infinite, maximal RA-suspensions (Euclidean $4$-circuits) splits the Coxeter suborbifold dual to $D$ along Euclidean suborbifolds into Seifert fibered pieces.  See Example \ref{ex:seifert}.

\begin{figure}[placement]
	\hspace{3.5cm}
	\input{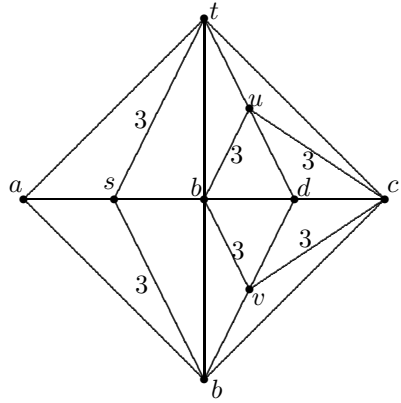}
	\caption{A Seifert disk in $L$}
	\label{fig:seifertdisk}
\end{figure}

\begin{Example}\label{ex:seifert} Suppose $L$ contains the Seifert subcomplex $D$ pictured in Figure \ref{fig:seifertdisk}, where the labels not denoted are $2$.  $D$ is made up of five maximal RA-suspensions: the suspension of $a$, $b$ and $c$, the suspension of $t$, $d$, and $b$, and the stars of $s$, $t$ and $v$.  Each of these intersect another along a Euclidean $4$-circuit.  The suborbifolds corresponding to these $4$-circuits have an intermediate cover isomorphic to a collared torus (an atoroidal piece).  The result is a splitting of the corresponding Coxeter suborbifold where each piece is either Seifert fibered or atoroidal.  
\end{Example}

\subsection{Applying Andreev's Theorem}\label{ss:andreev}
Since $L$ is metric flag, no two $3$-Euclidean vertices are connected by an edge.  So, stars of $3$-Euclidean vertices do not intersect in a $2$-simplex.

Let $s$ be a $3$-Euclidean vertex and let $T$ be an infinite, maximal RA-suspension.  Since $T$ is infinite, $\St_L(s)$ intersecting $T$ in one $2$-simplex with $s$ a suspension point implies $\St_L(s)\subset T$.  Then $T=L$, a suspension of a $3$-gon.  $T$ does not intersect $\St_L(s)$ in two $2$-simplices, since two edges of $L_s$ cannot be labeled $2$.  $T$ does not intersect $\St_L(s)$ in one $2$-simplex with a suspension point in $L_s$, since $L$ triangulates $\BS^2$.  By the same reason, $\St_L(s)$ and $T$ cannot intersect along an interior edge of either one.  Thus, if $L$ is not the suspension of a $3$-gon, the star of a $3$-Euclidean vertex and a Seifert subcomplex may only intersect along their boundary edges.    

\begin{figure}[h]
	\hspace{3.5cm}
	\input{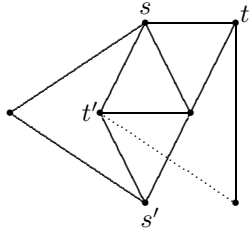}
	\caption{Two infinite RA-suspensions intersecting in a $2$-simplex.}
	\label{fig:two2-simplices1}
\end{figure}

Now let $S$ and $T$ denote infinite, maximal RA-suspensions of $L$.  Suppose $S$ and $T$ intersect in a $2$-simplex, sharing a suspension point.  Then since $L$ triangulates $\BS^2$, $S$ and $T$ intersect in an entire suspended edge and $S\cup T$ is an infinite RA-suspension, contradicting maximality.  Next, suppose $S$ and $T$ intersect in a $2$-simplex, but do not share a suspension point.  Since both $S$ and $T$ are infinite, both $S$ and $T$ suspend another vertex.  But since $L$ triangulates $\BS^2$, these vertices must coincide.  (See Figure \ref{fig:two2-simplices1}, where $\{s,s'\}$ and $\{t,t'\}$ denote the suspension points of $S$ and $T$ respectively.)  Then $S\cup T$ is an infinite RA-suspension where  $\{s,s',t,t'\}$ are the suspended points, again contradicting maximality.  Also because $L$ triangulates $\BS^2$, $S$ and $T$ cannot intersect along an interior edge of either one.  It follows that Seifert subcomplexes may only intersect along their boundary edges.  

\paragraph{Eliminating Seifert subcomplexes and the stars of $3$-Euclidean vertices.}  Suppose that $L$ is not the suspension of a $3$-gon.  Then from $L$, cut out each empty Euclidean $3$-circuit, each RA-cone on a Euclidean $3$-circuit and each Seifert subcomplex, capping off each remaining piece with a triangular cell or square cell.  (This operation is well-defined by the previous paragraphs.)  Let $\left[L-T\right]$ denote one of the remaining components.  (Here, $L-T$ denotes the corresponding full subcomplex of $L$.)  Then $\left[L-T\right]$ is a well-defined $2$-dimensional cell complex homeomorphic to $\BS^2$ with triangular and square faces in the strict sense that any non-empty intersection of two cells is a cell.  

The fact that $\left[L-T\right]$ is homeomorphic to $\BS^2$ is clear.  To see that any non-empty intersection of two cells is a cell, we need only check the intersection of the added triangular and square cells.  Two added squares cannot share all four edges, since then RA-suspensions in each removed Seifert subcomplex would intersect in a Euclidean $4$-circuit.  $L$ is not the suspension of a $3$-gon, so added triangles cannot share all their edges.  $\left[L-T\right]$ cannot be a 'pouch' (two vertices of $L$ would be connected by two different edges of $L$), so three edges of added squares cannot coincide, nor can two edges of added triangles.  By the same reasoning, two opposite edges of added squares cannot coincide in the manner of a cylinder.  $L$ is homeomorphic to $\BS^2$, so opposite edges of two added squares cannot intersect in the manner of a M\"obius band.  Two edges of an added triangle and two adjacent edges of added squares cannot coincide since no Euclidean $3$-circuit can have two of its edges labeled $2$.  Lastly, suppose adjacent edges of two added squares coincide.  Then the RA-suspension of three points, the point at the intersection of these adjacent edges and the opposite corners of the two squares, combines the two removed Seifert subcomplexes into a single Seifert subcomplex.  We have the following theorem.  

\begin{Theorem}\label{t:do2:10.3.1} Suppose that $L$ is not the boundary of a $3$-simplex and not a suspension of a $3$-gon.  Let $\left[L-T\right]$ be as above.  Then $K_{\left[L-T\right]}$ can be realized as a (possibly ideal), convex polytope in $\BH^{3}$.  (The ideal vertices correspond to the square or added triangular faces of $\left[L-T\right]$.)
  
$K_{\left[L-T\right]}$ can then be realized as a suborbifold of $K_L$ by cutting off the ideal vertices and labeling the resulting faces with the trivial group.    
\end{Theorem}

\begin{proof} If $K_{\left[L-T\right]}$ is a $3$-simplex, then $W_{L-T}$ is one of the non-compact hyperbolic Coxeter groups shown in Figure $6.3$ of \cite{humphreys}.  It acts on $\BH^3$ with fundamental chamber $K_{\left[L-T\right]}$, a finite volume simplex with ideal vertices dual to the added triangular faces of $\left[L-T\right]$.  Otherwise, we prove that $K_{\left[L-T\right]}$ satisfies the conditions of Andreev's theorem.

$\left[L-T\right]$ is a cell-complex with triangular and square faces, so $K_{\left[L-T\right]}$ has no more than three or four faces meeting at any vertex.  Condition (\ref{i:i}) is immediate under our hypothesis.  The remaining conditions refer to certain configurations of faces of the polytope.

$\left[L-T\right]$ contains no Euclidean $3$- or $4$-circuits not the boundary of a single face or the boundary of adjacent simplices , so it follows that every $3$- or $4$-prismatic element in $K_{\left[L-T\right]}$ satisfies condition (\ref{i:ii}) or (\ref{i:iii}).  $L$ is not the suspension of a $3$-gon, so the only way $K_{\left[L-T\right]}$ can be a triangular prism is if $\left[L-T\right]$ is the suspension of a $3$-gon where some face is an added triangular cell.  The boundary of this added triangle has edge labels $m_1$, $m_2$ and $m_3$ with the property that $\frac{\pi}{m_1}+\frac{\pi}{m_2}+\frac{\pi}{m_3}=\pi$.  So we know that not every suspension line is labeled $2$.  Thus, $K_{\left[L-T\right]}$ satisfies condition (\ref{i:iv}).  

To verify condition (\ref{i:v}) we note that if two faces $F_{1}$ and $F_{3}$ of $K_{\left[L-T\right]}$ intersect at a vertex, but are not adjacent, then this vertex must have valence $4$.  So this vertex corresponds to a square cell of $\left[L-T\right]$, where each edge is labeled $2$, and the two faces are dual to opposite corners $f_{1}$ and $f_{3}$ of the square.  The configuration in condition (\ref{i:v}) has a third face, $F_{2}$, adjacent to both the previous two.  So its dual vertex, $f_{2}$, is connected to both $f_1$ and $f_3$ in $\left[L-T\right]$.  Suppose both $m_{f_1f_2}$ and $m_{f_2f_3}$ equal $2$.  Then the boundary of the square, along with these edges form an infinite RA-suspension intersecting a removed Seifert subcomplex in a Euclidean $4$-circuit; a contradiction.  
\end{proof}

%

\subsection{The main result}
We are now ready to prove the main result, analogous to the JSJ-decomposition for $3$-dimensional manifolds (\cite{jash,joh}).  But first, we handle the cases $L$ is the boundary of a $3$-simplex or the suspension of a $3$-gon.  

\paragraph{The case where $L$ is the boundary of a $3$-simplex.} If $L$ is the boundary of a $3$-simplex, then $K_{L}$ is a $3$-simplex and we are unable to apply Andreev's theorem.  However, one can check that in this case $W_{L}$ is one of the groups listed in Figure $2.2$ or $6.2$ of \cite{humphreys} (n=4) and in fact, $\gS_L=\BE^3$ or $\gS_L=\BH^3$.  Thus, the characteristic suborbifold is either all of $K_L$ or it is empty.

\paragraph{$L$ the suspension of a $3$-gon.}  If $L$ is the suspension of a $3$-gon, the only conditions $K_{L}$ may fail to meet are (\ref{i:ii}) and (\ref{i:iv}).  If both (\ref{i:ii}) and (\ref{i:iv}) fail, then $\gS_L=\BR^3$.  If only (\ref{i:iv}) fails, then $L$ is a RA-suspension and $\gS_L=\BH^2\times\BE$.  In both of these cases the characteristic suborbifold is all of $K_L$.  If only (\ref{i:ii}) fails then $L$ decomposes as the union of two stars of $3$-Euclidean vertices.  If one of these stars is a RA-cone, then the characteristic suborbifold of $K_L$ is defined by this star.  If neither star is a RA-cone, the characteristic suborbifold is defined by the Euclidean $3$-circuit.  If both (\ref{i:ii}) and (\ref{i:iv}) are satisfied, $K_L$ is hyperbolic and the characteristic suborbifold is empty. 
  
\begin{Mainthm}\label{t:char} Let $K_L$ be a closed, irreducible, $3$-dimensional Coxeter orbifold.  There is a unique minimal collection of Euclidean suborbifolds of $K_L$ forming the boundary of the characteristic suborbifold of $K_L$, whose complement is atoroidal.  
\end{Mainthm}
\begin{proof} We have already described the characteristic suborbifold in the cases $L$ is the suspension of a $3$-gon or the boundary of a $3$-simplex.  In all other cases, cut out from $L$ each empty Euclidean $3$-circuit, each RA-cone on a Euclidean $3$-circuit and each Seifert subcomplex, as in Section \ref{ss:andreev}.  In $K_L$, this corresponds to cutting along triangular and quadrangular prismatic elements (three or four faces intersecting cyclically) whose cross-sections are incompressible $2$-dimensional Euclidean suborbifolds.  

The characteristic suborbifold is comprised of the suborbifolds of $K_L$ defined by the empty Euclidean $3$-circuits, the RA-cones on Euclidean $3$-circuits, and the Seifert subcomplexes.  Their boundary components are $2$-dimensional, Euclidean suborbifolds.  Since for a given $L$ these subcomplexes are uniquely defined and since Seifert subcomplexes are maximal, this set of boundary components is unique and minimal.  By Theorem \ref{t:do2:10.3.1}, the complement of the characteristic suborbifold is atoroidal.  
\end{proof}

As described in Section \ref{s:geom}, Seifert subcomplexes split along Euclidean $4$-circuits into infinite maximal RA-suspensions.  Thus, we have the following corollary to Theorem \ref{t:char}, geometrizing a Coxeter orbifold.  
\begin{Cor}\label{c:geom} Every closed, irreducible, $3$-dimensional Coxeter orbifold $K_L$ has a canonical decomposition along Euclidean suborbifolds into pieces which have the geometric structure of $\BH^3$, $\BE^3$, $\BE^2\times\left[-1,1\right]$ or $\BH^2\times \BE^1$.    
\end{Cor}
%

\section{Singer's conjecture}
Let $L$ be a metric flag simplicial complex, and let $A$ be a full subcomplex of $L$.  The following notation will be used throughout the remainder of the paper.
\begin{align}
\mfh_i(L) &:= \cH_i(\gS_L)\label{e:not1}\\
\mfh_i(A) &:= \cH_i(W_L\gS_A)\label{e:not2}\\
\gb_{i}(A)&:= \dim_{W_L}(\mfh_i(A)).\label{e:not4}
\end{align}
Here $\dim_{W_L}(\mfh_i(A))$ is the von Neumann dimension of the Hilbert $W_L$-module $W_L\gS_A$ and $\gb_{i}(A)$ is the $i^{\text{th}}$ $\ltwo$-Betti number of $W_L\gS_A$.  The notation in \ref{e:not2} and \ref{e:not4} will not lead to confusion since $\dim_{W_L}(W_L\gS_A)=\dim_{W_A}(\gS_A)$.  (See \cite{do2} and \cite{eckmann}).  We say that $A$ is \emph{$\ltwo$-acyclic}, if $\gb_{i}(A)=0$ for all $i$.

\subsection{Previous results in $\ell^2$-homology}\label{ss:previous}
\paragraph{Bounded geometry.}  The following result is proved by Cheeger and Gromov in \cite{cheeggrom}.  Suppose that $X$ is a complete contractible Riemannian manifold with uniformly bounded geometry (i.e. its sectional curvature is bounded and its injectivity radius is bounded away from $0$.)  Let $\Gamma$ be a discrete group of isometries on $X$ with $\Vol(X/\Gamma)<\infty$.  Then $\dim_{\Gamma}(\cH_k(\underline{E}\Gamma))=\dim_{\Gamma}(\cH_{k}(X))$, where $\underline{E}\Gamma$ denotes a universal space for proper $\Gamma$ actions, and $\cH_{k}(X)$ denotes the space of $\Ltwo$-harmonic forms on $X$.  Of particular interest to us is the case where $X=\BH^{3}$.  For it is proved by Dodziuk in \cite{dodz} that the $\Ltwo$-homology of any odd-dimensional hyperbolic space vanishes.   

\paragraph{Euclidean Space.}  The Cheeger Gromov result also implies that if $\gS_{L}=\BR^n$ for some $n$, then $\mfh_{\ast}(L)$ vanishes.  

\paragraph{Joins and suspensions.}  If a full subcomplex $A$ is the join of $A_1$ and $A_2$, i.e. $A=A_1\ast A_2$, where each edge connecting a vertex of $A_1$ with a vertex of $A_2$ is labeled $2$, then $W_A=W_{A_1}\times W_{A_2}$ and $\gS_A=\gS_{A_1}\times\gS_{A_2}$.  We may then use K\"unneth formula to calculate the $\ltwo$-Betti numbers of $\gS_{A}$; i.e.   
\begin{equation}\label{e:rt-angledjoin}
	\gb_{k}(A_{1}\ast A_{2})=\sum_{i+j=k}\gb_{i}(A_{1})\gb_{j}(A_{2}).
\end{equation}

\subsection{Singer's Conjecture for Coxeter groups}
Using the geometrical results proved in Section \ref{s:geom} and the results in Section \ref{ss:previous}, we can calculate the $\ltwo$-homology of a Coxeter system whose nerve $L$ is a triangulation of $\BS^2$.

\begin{Lemma}\label{l:euclidean} Suppose $L$ is not the boundary of a $3$-simplex and let $T\subset L$ denote an empty Euclidean $3$-circuit, the RA-cone on a $3$-Euclidean vertex, or an infinite maximal RA-suspension.  Then $T$ is $\ltwo$-acyclic.    
\end{Lemma}
\begin{proof} $T$ is full in $L$ by Lemma \ref{l:full} or \ref{l:rasfull}.  If $T$ is an empty Euclidean $3$-circuit, then $\gS_T=\BE^2$.  If $T$ is the RA-cone on a $3$-Euclidean vertex, then $T$ has the geometry of $\BE^2\times\left[-1,1\right]$.  In these cases, the result follows from the fact that the reduced $\ltwo$-homology of Euclidean space vanishes (Section \ref{ss:previous}) and equation (\ref{e:rt-angledjoin}).  If $T$ is a RA-suspension, where $P$ denotes the suspension points, then the result follows from equation (\ref{e:rt-angledjoin}), since $\gS_P=\BE$ and thus $\gb_{\ast}(P)=0$.  
\end{proof}

As a result, Singer's Conjecture for Coxeter groups in dimension $3$ follows from Corollary \ref{c:geom}.  
\begin{Cor}\label{t:singer} Let $L$ be a metric flag triangulation of $\BS^2$.  Then $\mfh_{i}(L)=0$ for all $i$.
\end{Cor}
\begin{proof} If $L$ is the boundary of a $3$-simplex, then $\gS_L$ is either $\BE^3$ or $\BH^3$.  In either case, $L$ is $\ltwo$-acyclic (Section \ref{ss:previous}).  Otherwise, decompose $K_L$ into geometric pieces as in Corollary \ref{c:geom}.  The subcomplexes defining the non-hyperbolic pieces are $\ltwo$-acyclic by Lemma \ref{l:euclidean}.  Those defining the hyperbolic pieces are $\ltwo$-acyclic by Theorem \ref{t:do2:10.3.1} and the results in Section \ref{ss:previous}.  Since the intersection of these subcomplexes define Euclidean suborbifolds, the result follows from Mayer-Vietoris.  
  
%
\end{proof}

\bibliography{mybib}   

\end{document}